\documentclass[12pt]{amsart}

\newtheorem{theorem}{Theorem}

\newtheorem{thm}[theorem]{Theorem}

\newtheorem{proposition}[theorem]{Proposition}

\theoremstyle{definition}

\newtheorem{df}[theorem]{Definition}

\newtheorem{problem}[theorem]{Problem}

\theoremstyle{remark}

\newcommand{\alg}{\mathbf}
\def\v{\mathcal V}  
\DeclareMathOperator{\con}{Con}
\newcommand{\cona}{\con {\mathbf A}}

\title[\lowercase{$m$}-permutable varieties, congruence identity]
{A congruence identity satisfied by $m$-permutable varieties}

\author[]{Paolo Lipparini} 
\address{Dipartimento di Matematica\\
Viale della Ricerca Scientifica\\
II Universit\`a di Rom (Tor Vergata)\\
I-00133 ROME 
ITALY
}
\email{lipparin@axp.mat.uniroma2.it}
\urladdr{http://www.mat.uniroma2.it/\textasciitilde lipparin}
\thanks{The author has received support from MPI and GNSAGA.
We wish to express our  gratitude to K. Kearnes for stimulating
correspondence} 


\begin{document} 

\begin{abstract} 
We present a new and useful
 congruence identity satisfied by $m$-permutable varieties.
\end{abstract}

\maketitle 


It has been proved in \cite{Lp1}  that every $m$-permutable variety satisfies a
non-trivial lattice identity (depending only on $m$).
In \cite{lpabn} we have found another interesting identity:

\begin{thm}\label{npermabn}
For $m\geq 3$, every $m$-permutable variety 
satisfies the congruence identity
$ \alpha \beta _h = \alpha \gamma _h$,
for $h= m[\frac{m+1}{2}]-1$   
\end{thm} 
 
Here, $[\ ]$ denotes integer part, and $\beta _h$, $\gamma _h$
are defined as usual:
\[
 \beta_0= \gamma_0= 0,
\]
\[
\beta_{n+1}=\beta+\alpha\gamma_n, \qquad \qquad 
\gamma _{n+1}=\gamma +\alpha\beta _n.
\]

The proof of Theorem \ref{npermabn}
consists of two steps. 
As for the second step, it is an easy application of commutator theory,
but in the present note we shall be concerned only with
the first step.

The first step
is a commutator-free proof of the following

\begin{thm}\label{npermimplXasquare}
If every subalgebra of $\alg A^2$ is $m$-permutable 
then $\alg A$ satisfies $(X_m)$.
More generally:

(i) If every congruence of $ {\alg A} $, thought of as a
 subalgebra of $\alg A^2$, is $m$-permutable 
then $\alg A$ satisfies $(X_m)$.

(ii)  If every subalgebra of $\alg A^2$ 
generated by $m+1$ elements
is $m$-permutable 
then $\alg A$ satisfies $(X_m)$;
actually, ${\alg A} $ satisfies the stronger version
of $(X_m)$ in which $\alpha $ is only supposed to be a compatible
relation on ${\alg A} $, and $ \delta $ is any relation on ${\alg A} $.  
\end{thm} 

In the statement of Theorem  \ref{npermimplXasquare}
we have used:

\begin{df}\label{Xm} 
If $\alpha,\beta,\gamma, \delta  $ are congruences on some algebra,
and $m$ is a natural number, we shall denote by $(X_m)$  
the following identity.
\begin{gather*}
\alpha ( \beta \circ \alpha ( \gamma \circ \alpha ( \beta \circ \dots
\alpha ( \gamma ^\bullet \circ \alpha ( \beta ^\bullet \circ \alpha \delta  \circ \beta ^\bullet) \circ \gamma ^\bullet) \ldots
 \circ \beta )\circ \gamma) \circ \beta ) = 
 \\ 
\alpha ( \gamma  \circ \alpha ( \beta  \circ \alpha ( \gamma  \circ \dots
\alpha ( \beta  ^\bullet \circ \alpha ( \gamma  ^\bullet \circ \alpha \delta  \circ \gamma  ^\bullet) \circ \beta  ^\bullet) \ldots
 \circ \gamma  )\circ \beta ) \circ \gamma  ) \ 
\end{gather*}
with exactly $m$ open brackets
 (and exactly $m$ closed brackets) on each side,
and where $ \beta ^\bullet= \beta $, $ \gamma ^\bullet= \gamma $  if 
$m$ is odd, and   
$ \beta ^\bullet= \gamma  $, $ \gamma ^\bullet= \beta  $  if 
$m$ is even.
\end{df}

Thus, $(a_0,b_0)$ belongs to the left-hand side of $(X_m)$
if and only if  there are further elements 
$ a_1, a_2, \dots, a_m$ and 
$ b_1, b_2, \dots, b_m$
such that 
\begin{align*}
&& &a_i \alpha b_i, &&\text{for  }i=0,\dots,m, \\
&& &a_m \delta b_m, \\
&a_i  \beta a_{i+1},  & &b_i  \beta b_{i+1}, && \text{for } i \text{ even, } 0\leq i\leq m-1, \\
&a_i \gamma a_{i+1},  & &b_i \gamma b_{i+1}, && \text{for } i \text{ odd, } 0\leq i\leq m-1. 
\end {align*} 
The conditions asserting that $(a_0,b_0)$ belongs to the right-hand side of $(X_m)$
are similar, with $\beta $ and $ \gamma $ interchanged.
See \cite{lpabn} for more comments about Definition
\ref{Xm}, and for a picture.

Our original proof of Theorem \ref{npermabn}
used Hagemann and Mitschke's
terms \cite{HM} and in the first step we only got
the following weaker version of Theorem \ref{npermimplXasquare}: 

\begin{thm}\label{npermimplX} 
Every $m$-permutable variety satisfies $(X_m)$. 
\end{thm} 

Then we found a simpler proof of the stronger
Theorem \ref{npermimplXasquare}.
However, our original argument using
Hagemann and Mitschke's
terms provides the following variant, which does not follow
from Theorem \ref{npermimplXasquare},
and which we prove here in the hope for further applications.
There is even a more general version which we do not state
here.

\begin{proposition}\label{XX} 
Every algebra $\alg A$ belonging to an $m$-permutable variety satisfies 
\begin{multline*}
R_{0} ( S_{1} \circ R_{1} ( S_{2} \circ R_{2} ( S_{3} \circ \dots
R_{m-2} ( S_{m-1}  \circ R_{m-1} ( S_{m}  \circ R_{m}  \circ T_{m} ) \circ T_{m-1} ) \dots \\
 \circ T_{3} )\circ T_{2}) \circ T_{1} ) \ \subseteq \  
R
_{0} ( S'_{1} \circ R'_{1} ( S'_{2} \circ R'_{2} ( S'_{3} \circ \dots \\
R'_{m-2} ( S'_{m-1}  \circ R'_{m-1} ( S'_{m}  \circ R
_{m}  \circ T'_{m} ) \circ T'_{m-1} ) \ldots
 \circ T '_{3} )\circ T'_{2}) \circ T'_{1} )
\end{multline*}
for every relations
$ R_{0}, \dots, R_{m}$, 
and reflexive relations
$S_{1}, \dots, S_{m}, T_{1},
$ $ 
 \dots, T_{m}$  
on ${\alg A}$, where
\begin{gather*}
R'_{i} = \overline{R_{i-1} \cup  R_{i} \cup  R_{i+1}},
\qquad \text{ for } i=1, \dots, m-1,\\
S'_{1} = \overline {S_{2}}, \qquad 
T'_{1} = \overline {T_{2}}, \qquad
S'_{m} =  \overline {S_{m-1}}, \qquad 
T'_{m} =  \overline {T_{m-1}},\\
S'_{i} = \overline{S_{i-1}} \circ \overline{S_{i+1}},
 \quad T'_{i} = \overline{T_{i+1}} \circ \overline{T_{i-1}},
\quad \text{ for } i=2, \dots, m-1,
\end{gather*}
and $\overline{X}$ denotes the least compatible relation
containing $X$.
\end{proposition}

Here are the notations we use:
$ \alpha, \beta   $ denote \textit{congruences}
on some \textit{algebra} ${\alg A}$. 
Join and meet in the lattice $ \cona$ of all congruences of 
${\alg A}$ are denoted, respectively, by $+$ and 
juxtaposition.
We use juxtaposition also to denote intersection. 

Relational product is denoted by $\circ$, 
and $ \alpha \circ_n \beta $ is a shorthand     
for $ \alpha \circ \beta \circ \alpha \circ \beta \circ \dots$,
with $n-1$ occurrences of $\circ$.   
Two congruences  $ \alpha,  \beta $
are said to {\em $m$-permute}
if and only if $ \alpha \circ_m \beta =\beta \circ_m \alpha$
(thus, in particular, $  \alpha + \beta =\alpha \circ_m \beta$).
An algebra ${\alg A}$ is {\em $m$-permutable}
if and only if every pair of congruences in 
${\alg A}$ $m$-permute.
A variety $\v$ is {\em $m$-permutable}
if and only if every algebra in  
$\v$ is $m$-permutable.

\begin{theorem}\label{hm}
\cite{HM}  A variety is $m$-permutable if and only if  it has terms
$t_0,\dots,t_{m}$ satisfying:

(i) $x=t_0(x,y,z)$,
  
(ii) $t_i(x,x,y)=t_{i+1}(x,y,y)$,  for $i=0,\dots,m-1$,

(iii) $t_{m}(x,y,z)=z$.
\end{theorem}

\begin{proof}[Proof of Theorem \ref{npermimplX}]
Let ${\alg A}$ belong to some
$m$-permutable variety, and $ \alpha,\beta,\gamma,\delta \in \cona$.

Suppose that $(a_0,b_0)$ belongs to the left-hand side
of $(X_m)$. 
It is enough to show that $(a_0,b_0)$
belongs to the right-hand side.
The reverse inclusion is obtained by  symmetry.

Suppose that there are elements $ a_1, a_2, \dots, a_m$ and 
$ b_1, b_2, \dots, b_m$
as in the comment after Definition \ref{Xm}.
We want to obtain elements 
 $ c_1, c_2, \dots, 
$ $ 
c_m$ and 
$ d_1, d_2, \dots, d_m$ witnessing that
$(a_0,b_0)$
belongs to the right-hand side.

We apply the usual arguments showing that
the terms from Theorem \ref{hm} imply $m$-permutability:
we apply the terms $t_0,\dots,t_{m}$ symmetrically
to the $a_i$'s and the $b_i$'s, thus getting everything congruent
modulo $ \alpha $.

Namely, we claim that the elements
$c_i=t_i(a_{i-1}, a_i, a_{i+1})$, 
$d_i=t_i(b_{i-1}, b_i, b_{i+1})$,  for
$i=1,\dots,m-1$, and $c_m=a_m$, $d_m=b_m$
witness that $(a_0,b_0)$ belongs to the right-hand side of $(X_m)$,
where the terms $t_i$ are given by Theorem \ref{hm}.

Indeed, for $i=1,\dots,m-1$, $c_i=t_i(a_{i-1}, a_i, a_{i+1}) \alpha  
t_i(b_{i-1}, b_i, b_{i+1})=d_i$ is trivial, since 
$a_i\alpha b_i$ for  $i=0, \dots, m$,
and since $ \alpha $ is compatible. 
Moreover, $c_m=a_m \alpha \delta  b_m=d_m$.

By Conditions (i) and (ii) in Theorem \ref{hm}, and since $a_1 \gamma a_2$, 
we have that 
$ a_0= t_0(a_{0}, a_0, a_{1})= t_1(a_{0}, a_1, a_{1}) \gamma 
t_1(a_{0}, a_1, a_{2})=c_1$. Similarly, $b_0 \gamma d_1$. 

By Conditions (ii) and (iii) in Theorem \ref{hm}, say, for $m$ odd, 
$ c_{m-1} = t _{m-1} (a_{m-2}, a_{m-1}, a_{m})
\gamma  t _{m-1} (a_{m-1}, a_{m-1}, a_{m})=
t _{m} (a_{m-1}, a_{m}, a_{m})=a_m=c_m$,
since $a_{m-2} \gamma a_{m-1}$.
Similarly, $d_{m-1} \gamma d_m$.
The case $m$ even is similar: just replace $ \gamma $  by $ \beta $.

Now, let $1 \leq i \leq m-2$ and, say, $i$ odd. 
By Condition (ii) in Theorem \ref{hm}, and since $ a_{i-1} \beta  a_i $,
and $ a_{i+1} \beta  a _{i+2}  $, we get
$c_i=t_i(a_{i-1}, a_i, a_{i+1}) \beta 
$ $ 
t_i(a_i, a_i, a_{i+1})=
t_{i+1} (a_{i}, a_{i+1}, a_{i+1}) \beta
t_{i+1} (a_{i}, a_{i+1}, a_{i+2})=c _{i+1}$.
Similarly, $d_{i} \beta  d_{i+1}$.

In the case $i$ even we get  
$c_{i} \gamma   c_{i+1}$ and 
$d_{i} \gamma   d_{i+1}$ in the same way.

In conclusion, we have showed that
the elements $c_i$, $d_i$ satisfy the desired relations.  
\end{proof} 

\begin{proof}[Proof of Proposition \ref{XX}]
The argument is essentially the same as above.
Suppose that
$(a_0,b_0)$ belongs to the left-hand side
of the inclusion in Proposition \ref{XX}.
This is witnessed by  elements 
$ a_1, a_2, \dots, a_m$ and 
$ b_1, b_2, \dots, b_m$
such that 
\begin{align*}
&& &a_i R_i b_i &&\text{for  }i=0,\dots,m \\
&a_i  S_{i+1} a_{i+1}  & &b_{i+1} T_{i+1} b_{i} && \text{for }  
i=0,\dots, m-1
\end {align*} 

We are going to show that the elements
$c_i=t_i(a_{i-1}, a_i, a_{i+1})$, 
$d_i=t_i(b_{i-1}, b_i, b_{i+1})$,  for
$i=1,\dots,m-1$, and $c_m=a_m$, $d_m=b_m$
witness that $(a_0,b_0)$ belongs to the right-hand side of 
the inclusion in Proposition \ref{XX},
where the terms $t_i$ are given by Theorem \ref{hm}.

Indeed, for $i=1,\dots,m-1$, 
$
c_i=
t_i(a_{i-1}, a_i, a_{i+1})
\overline{R_{i-1} \cup  R_{i} \cup  R_{i+1}}
$ $ 
t_i(b_{i-1},
b_i, b_{i+1})
=d_i
$
since 
$ a_{i-1} R_{i-1} b_{i-1} $, 
$ a_i R_i b_i $,
$ a_{i+1} R_{i+1} b_{i+1} $, and, by definition, 
$R'_i=\overline{R_{i-1} \cup  R_{i} \cup  R_{i+1}}$
is compatible and contains 
$R_{i-1}$, 
$ R_i $ as well as
$ R_{i+1}$.

Moreover, 
$a_0 R_0 B_0$, and  
$c_m=a_m R_m  b_m=d_m$.

By Conditions (i) and (ii) in Theorem \ref{hm}, and since 
$S$ is reflexive and
$a_1 S_2 a_2$, 
we have that 
$ a_0= t_0(a_{0}, a_0, a_{1})= t_1(a_{0}, a_1, a_{1}) 
\overline{S}_2 
t_1(a_{0}, a_1, a_{2})=c_1$. 
Thus, $a_0 S'_1 c_1$.
Symmetrically, 
$d_1= t_1(b_{0}, b_1, b_{2})
\overline{T}_2 
t_1(b_{0}, b_1, b_{1})=
t_0(b_{0}, b_0, b_{1})=b_0$.
Hence, $d_1 T'_1 b_0$. 

By Conditions (ii) and (iii) in Theorem \ref{hm},
$ c_{m-1} = t _{m-1} (a_{m-2}, a_{m-1}, 
$ $ 
a_{m}) \overline{S_{m-1}}
  t _{m-1} (a_{m-1}, a_{m-1}, a_{m})=
t _{m} (a_{m-1}, a_{m}, a_{m})=a_m=c_m$,
since $a_{m-2} S_{m-1} a_{m-1}$
and $S_{m-1}$ is reflexive.
Thus, $ c_{m-1} S'_m c_{m}$.
Symmetrically, $d_{m} T'_{m} d_{m-1}$.

Now, let $1 \leq i \leq m-2$. 
By Condition (ii) in Theorem \ref{hm}, and since $ a_{i-1} S_i  a_i $,
and $ a_{i+1} S_{i+2}  a _{i+2}  $, we get
$c_i=t_i(a_{i-1}, a_i, a_{i+1}) \overline{S_{i}} 
t_i(a_i, a_i, a_{i+1})=
t_{i+1} (a_{i}, a_{i+1}, a_{i+1}) \overline{S_{i+2}}
t_{i+1} (a_{i}, a_{i+1}, a_{i+2})=c _{i+1}$.
Thus, 
$c_i (\overline{S_{i}} 
\circ \overline{S_{i+2}})
c _{i+1}$, that is, 
$c_i S'_{i+1}c _{i+1}$,
by the definition of 
$S'_{i+1}$.
Symmetrically, $d_{i+1} T'_{i+1} d_{i}$.

In conclusion, we have showed that
the elements $c_i$, $d_i$ satisfy the desired relations.  
\end{proof}

Notice that Theorem \ref{npermimplX} is the particular case
$R_0=R_1=\dots=R_{m-1}= \alpha $, $R_m=\alpha\delta $,
$S_1=S_3=S_5=\dots =\beta $, $T_1=T_3=T_5= \dots =\beta $,
$S_2=S_4=S_6=\dots =\gamma  $, $T_2=T_4=T_6= \dots =\gamma  $
of Proposition \ref{XX}.

Our proof of Theorem
\ref{npermabn} gives something more.
Consider the following property $(X_m)^*$
which is weaker than $(X_m)$.
\begin{gather*}
\alpha ( \beta \circ \alpha ( \gamma \circ \alpha ( \beta \circ \dots
\alpha ( \gamma ^\bullet \circ \alpha ( \beta ^\bullet \circ \alpha \delta  \circ \beta ^\bullet) \circ \gamma ^\bullet) \ldots
 \circ \beta )\circ \gamma) \circ \beta ) \subseteq
 \\ 
\Big(
\alpha ( \gamma  \circ \alpha ( \beta  \circ \alpha ( \gamma  \circ \dots
\alpha ( \beta  ^\bullet \circ \alpha ( \gamma  ^\bullet \circ \alpha \delta  \circ \gamma  ^\bullet) \circ \beta  ^\bullet) \ldots
 \circ \gamma  )\circ \beta ) \circ \gamma  ) 
\Big)^*
\end{gather*}
with $ m$ normal-sized open parenthesis on each side, where $^*$ denotes transitive closure.


\begin{theorem}\label{xastabh}
Suppose that $\v$ is a variety satisfying 
$(X_m)^*$ 
(or just $(Y_{m-1})^*$, see below)
for some $m$.
In addition, suppose that on all
congruence lattices of algebras in $\v$ it is possible to define
a monotone and submultiplicative commutator operation $[-,-]$
such that for some $k>0$ the following properties hold:

(i) 
$[ \beta + \gamma  , \alpha ] \leq \alpha \beta _{k} $, 
for every algebra $ {\alg A} \in \v $ and congruences
$\alpha, \beta,\gamma  \in \cona$; and

(ii) there is a ternary term $d$ 
such that 
\[
d(b,b,a) \equiv a \equiv d(a,b,b) \pmod{ [  \alpha , \alpha ] }
\] 
for every algebra $ {\alg A} \in \v $, every congruence
$\alpha \in \cona$ and elements $a \alpha b \in A$.

Then $\v$ satisfies   
$ \alpha \beta _h = \alpha \gamma _h$,
for $h= k+m-1$.   
\end{theorem} 

We have a long technical proof
showing that if a variety satisfies $(X_m)^*$
then 
for some $n$ and $k$ the commutator
$[ \alpha ,\beta|n]$ satisfies
conditions (i) and (ii) in Theorem 
\ref{xastabh}. Thus we get:
{\em If a variety $ \v $ satisfies 
$(X_m)^*$
for some $m$ 
then there is some $h$ 
(depending on $\v$) such that 
$ \v $ satisfies 
$ \alpha \beta _h = \alpha \gamma _h$}.
 
Is the converse true?

\begin{problem}\label{xmiffabh} 
Is it true that if $ \v $ satisfies 
$ \alpha \beta _h = \alpha \gamma _h$ for some $h$
then 
$ \v $ satisfies 
$(X_m)^*$
for some $m$?
\end{problem} 

Notice that 
$(X_m)^*$
is equivalent to
\begin{gather*}
\Big(
\alpha ( \beta \circ \alpha ( \gamma \circ \alpha ( \beta \circ \dots
\alpha ( \gamma ^\bullet \circ \alpha ( \beta ^\bullet \circ \alpha \delta  \circ \beta ^\bullet) \circ \gamma ^\bullet) \ldots
 \circ \beta )\circ \gamma) \circ \beta ) \Big)^*
=
 \\ 
\Big(
\alpha ( \gamma  \circ \alpha ( \beta  \circ \alpha ( \gamma  \circ \dots
\alpha ( \beta  ^\bullet \circ \alpha ( \gamma  ^\bullet \circ \alpha \delta  \circ \gamma  ^\bullet) \circ \beta  ^\bullet) \ldots
 \circ \gamma  )\circ \beta ) \circ \gamma  ) 
\Big)^*
\end{gather*}

Consider the following condition 
 $(Y_m)$ 
\begin{gather*}
\alpha ( \beta \circ \alpha ( \gamma \circ \alpha ( \beta \circ \dots
\alpha ( \gamma ^\bullet \circ \alpha ( \beta ^\bullet \circ \alpha \gamma ^\bullet  \circ \beta ^\bullet) \circ \gamma ^\bullet) \ldots
 \circ \beta )\circ \gamma) \circ \beta ) = 
 \\ 
\alpha ( \gamma  \circ \alpha ( \beta  \circ \alpha ( \gamma  \circ \dots
\alpha ( \beta  ^\bullet \circ \alpha ( \gamma  ^\bullet \circ \alpha \beta  ^\bullet  \circ \gamma  ^\bullet) \circ \beta  ^\bullet) \ldots
 \circ \gamma  )\circ \beta ) \circ \gamma  ) \ 
\end{gather*}
with exactly $m$ open brackets
 on each side,
and where $ \beta ^\bullet= \beta $, $ \gamma ^\bullet= \gamma $  if 
$m$ is odd, and   
$ \beta ^\bullet= \gamma  $, $ \gamma ^\bullet= \beta  $  if 
$m$ is even.

Notice that $(X_m)$ implies $(Y_{m-1})$: just take $ \delta =0$. 
In fact, in our proof of Theorem \ref{npermabn}
it is enough to assume $(Y_{m-1})$ in place of $(X_m)$.

We can also consider $(Y_m)^*$
\begin{gather*}
\alpha ( \beta \circ \alpha ( \gamma \circ \alpha ( \beta \circ \dots
\alpha ( \gamma ^\bullet \circ \alpha ( \beta ^\bullet \circ \alpha \gamma ^\bullet  \circ \beta ^\bullet) \circ \gamma ^\bullet) \ldots
 \circ \beta )\circ \gamma) \circ \beta ) \subseteq
 \\ 
\Big(
\alpha ( \gamma  \circ \alpha ( \beta  \circ \alpha ( \gamma  \circ \dots
\alpha ( \beta  ^\bullet \circ \alpha ( \gamma  ^\bullet \circ \alpha \beta  ^\bullet  \circ \gamma  ^\bullet) \circ \beta  ^\bullet) \ldots
 \circ \gamma  )\circ \beta ) \circ \gamma  ) 
\Big)^*
\end{gather*}
with $ m$ normal-sized open parenthesis on each side.

$(Y_m)^*$ is equivalent to
\begin{gather*}
\Big(
\alpha ( \beta \circ \alpha ( \gamma \circ \alpha ( \beta \circ \dots
\alpha ( \gamma ^\bullet \circ \alpha ( \beta ^\bullet \circ \alpha \gamma ^\bullet  \circ \beta ^\bullet) \circ \gamma ^\bullet) \ldots
 \circ \beta )\circ \gamma) \circ \beta ) 
\Big)^*
=
 \\ 
\Big(
\alpha ( \gamma  \circ \alpha ( \beta  \circ \alpha ( \gamma  \circ \dots
\alpha ( \beta  ^\bullet \circ \alpha ( \gamma  ^\bullet \circ \alpha \beta  ^\bullet  \circ \gamma  ^\bullet) \circ \beta  ^\bullet) \ldots
 \circ \gamma  )\circ \beta ) \circ \gamma  ) 
\Big)^*
\end{gather*}
Notice that $(X_m)^*$ implies $(Y_{m-1})^*$.

\end{document}